\documentstyle[12pt]{article}

\leftmargin\leftmargini
\labelwidth\leftmargini\advance\labelwidth-\labelsep

\def\R{\relax\ifmmode I\!\!R\else$I\!\!R$\fi}

\def\Z{\relax\ifmmode Z\!\!\!Z\else$Z\!\!\!Z$\fi}

\def\C{\relax\ifmmode C\!\!\!\!I\else$C\!\!\!\!I$\fi}

\def\K{\relax\ifmmode I\!\!K\else$I\!\!K$\fi}

\def\N{\relax\ifmmode I\!\!N\else$I\!\!N$\fi}

\newcounter{defcounter}[section]
{\vspace{0.1cm}\begin{sloppypar}\noindent\stepcounter{defcounter}{\bfseries Definition
      \thesection.\thedefcounter}}%
{\end{sloppypar}\vspace{0.1cm}}

\newtheorem{theorem}{Theorem}[section]
\newtheorem{proposition}{Proposition}[section]

\newcommand{\proof}{{\noindent\bf Proof. }}

\begin{document}
\thispagestyle{empty}
\begin{center}
{\Large {\bf Dimension theoretical properties of generalized Baker's Transformations}}
\end{center}
\begin{center}
J. Neunh\"auserer\\
Institut f\"ur Theoretische Physik,\\ Technische Universit\"at Clausthal.
\footnote{This work is supported by the DFG research group 'Zeta functions
and locally symmetric spaces'.}\\
Arnold-Sommerfeld-Str. 6, \\ 38678 Clausthal-Zellerfeld, Germany\\
joerg.neunhaeuserer@tu-clausthal.de
\end{center}
\begin{abstract}
We show that for generalized Baker's transformations there is a
parameter domain where we have an absolutely continuous ergodic measure
and in direct neighborhood there is a parameter domain where not even the variational principle
for Hausdorff dimension holds.    \\
  {\bf MSC 2000: 37C45, 28A80, 28D20}
\end{abstract}

\section{Introduction \label{Chap1}}
In the modern theory of dynamical systems geometrical invariants
like Hausdorff and box-counting dimension of invariant sets and measures
seems to have their place beside classical invariants like entropy
and Lyapunov exponents.  In the last decades a dimension theory of
dynamical systems was developed and we have general results for
conformal systems, see \cite{[PE]}, \cite{[FA]}  and references therein. On the other hand the existence of different rates of contraction
or expansion in different directions forces mathematical problems that
are not completely solved. We have general results on hyperbolic measures (\cite{[LY]}, \cite{[BPS]}) but the question
if there exists an ergodic measure with full dimension (the dimension of a given invariant set) is only solved in special cases (\cite{[KP]}, \cite{[GL]}, \cite{[MC]}, \cite{[SCH]}).
In this work we consider a generalization of the Baker's transformation, a
simple example of a 'chaotic' dynamical system that may be found in many
standard text books \cite{[SU]}. In the case that the transformations are invertible dimensional
theoretical properties are fairly easy to understand and the results
seem to be folklore in the dimension theory of dynamical systems (see Theorem 2.1). We will be interested here in the case that the transformations are not invertible.  Our main result describes a phebomenon which was in this form not
observed before. In fact there is a parameter domain were there generically exists an absolutely continuous ergodic measure which obviously has full dimension on the attractor (see Theorem 2.2).
On the other hand in the neighborhood there ist a parameter domain were the
variational principle for Hausdorff dimension does not hold, the dimension
of the attractor can not even be approximated by the dimension of ergodic measures (see Theorem 2.2). A kind of bifurcation occurs. Also we illustrate
this phenomenon only in a very simple case we think it may generically occur
for endomorphisms.\\
The rest of the paper is organized as follows. In section two
we introduce the systems we study and present our main results. In section three we find a symbolic coding for the dynamics of generalized Baker's transformations through a factor of full shift on two symbols and represent all ergodic measures
using this coding. In section four we construct absolutely continues ergodic measure
for generalized Baker's transformation using our results on
overlapping self-similar measure \cite{[NE1]} and thus proof Theorem 2.2~. In section
we find upper estimates on the dimension of all ergodic measure and proof
Theorem 2.3~.

\section{Notations and results}
We define a {\it generalized Baker's Transformation} on the square by
\[          f_{\beta_{1},\beta_{2}} :           [-1,1]^2 \longmapsto  [-1,1]^2\]
\[ f_{\beta_{1},\beta_{2}}(x,y)=
\lbrace
\begin{array}{cc}
 (\beta_{1} x+(1-\beta_{1}),2y -1) \quad \mbox{if} \quad y \ge 0   \\
 (\beta_{2} x-(1-\beta_{2}),2y +1 )\quad \mbox{if}\quad y < 0
\end{array} \]
for parameter values $\beta_{1},\beta_{2}\in (0,1)$.
We call this family of maps generalized Baker's transformations
because if we set $\beta=\beta_{1}=\beta_{2}$ we get the
class of Baker's transformation studied by Alexander and York.
and for $\beta=0.5$
we get the well known
classical Baker's transformation \cite{[AY]}.\newline\newline
Let us first consider the case $\beta_{1}+\beta_{2}<1$
. In this case the attractor of the map $f_{\beta_{1},\beta_{2}}$
\[\Lambda_{\beta_{1},\beta_{2}}=\bigcap_{n=0}^{\infty} f^{n}_{\beta_{1},\beta_{2}}([-1,1]^2)\]
is a product of a cantor set with the interval $[-1,1]$ and the dimensional
theoretical properties of the system are easy to deduce.
\begin{theorem}
Let $\beta_{1}+\beta_{2}<1$ and $d$ be the unique positive number
satisfying \[ \beta_{1}^{d}+\beta_{2}^{d}=1 \] then
\[ \dim_{B} \Lambda_{\beta_{1},\beta_{2}}=\dim_{H} \Lambda_{\beta_{1},\beta_{2}}=d+1\] and there is an $f_{\beta_{1},\beta_{2}}$-ergodic measure $\mu$ of full dimension i.e. $\dim_{B}\mu=\dim_{H} \mu=d+1$.
\end{theorem}
This result seems to be folklore in the dimension theory of dynamical
systems. The box-counting dimension of $\Lambda_{\beta_{1},\beta_{2}}$
is easy to calculate and the ergodic measure of full dimension is
constructed as a product of a Cantor measure with weights $(\beta_{1}^{d},\beta_{2}^{d})$ on the real line with the normalized Lebesgue measure on $[-1,1]$. See section 23 of \cite{[PE]}) for this martial.
\begin{center}
\unitlength=1.00mm
\special{em:linewidth 0.4pt}
\linethickness{0.4pt}
\unitlength=1.00mm
\special{em:linewidth 0.4pt}
\linethickness{0.4pt}
\begin{picture}(150.00,78.00)
\put(10.00,8.00){\dashbox{2.00}(60.00,35.00)[cc]{}}
\put(10.00,43.00){\framebox(60.00,35.00)[cc]{}}
\put(90.00,8.00){\dashbox{1.41}(20.00,70.00)[cc]{}}
\put(116.00,8.00){\framebox(34.00,70.00)[cc]{}}
\put(72.00,43.00){\vector(1,0){16.00}}
\put(79.00,40.00){\makebox(0,0)[cc]{$f_{\beta_{1},\beta_{2}}$}}
\put(135.00,4.00){\makebox(0,0)[cc]{$2\beta_{1}$}}
\put(102.00,4.00){\makebox(0,0)[cc]{$2\beta_{2}$}}
\put(139.00,4.00){\vector(1,0){12.00}}
\put(130.00,4.00){\vector(-1,0){15.00}}
\put(105.00,4.00){\vector(1,0){5.00}}
\put(98.00,4.00){\vector(-1,0){7.00}}
\end{picture}
\end{center}
{\bf Figure 1:} The action of $f_{\beta_{1},\beta_{2}}$ on the square $[-1,1]^{2}$ in the case $\beta_{1}+\beta_{2}<1$\newline\newline
Now consider the case $\beta_{1}+\beta_{2}\ge1$ the attractor is obviously the whole square $[-1,1]^{2}$ which has
Hausdorff and box-counting dimension two.
\begin{center}
\unitlength=1.00mm
\special{em:linewidth 0.4pt}
\linethickness{0.4pt}
\unitlength=1.00mm
\special{em:linewidth 0.4pt}
\linethickness{0.4pt}
\begin{picture}(150.00,78.00)
\put(10.00,8.00){\dashbox{2.00}(60.00,35.00)[cc]{}}
\put(10.00,43.00){\framebox(60.00,35.00)[cc]{}}
\put(90.00,8.00){\dashbox{1.41}(45.00,70.00)[cc]{}}
\put(116.00,8.00){\framebox(34.00,70.00)[cc]{}}
\put(72.00,43.00){\vector(1,0){16.00}}
\put(79.00,40.00){\makebox(0,0)[cc]{$f_{\beta_{1},\beta_{2}}$}}
\put(135.00,4.00){\makebox(0,0)[cc]{$2\beta_{1}$}}
\put(110.00,1.00){\makebox(0,0)[cc]{$2\beta_{2}$}}
\put(139.00,4.00){\vector(1,0){12.00}}
\put(130.00,4.00){\vector(-1,0){15.00}}
\put(113.00,1.00){\vector(1,0){22.00}}
\put(106.00,1.00){\vector(-1,0){17.00}}
\end{picture}
\end{center}
{\bf Figure 2:} The action of $f_{\beta_{1},\beta_{2}}$ on the square $[-1,1]^{2}$ in the case $\beta_{1}+\beta_{2}>1$\newline\newline
The interesting problem in this situation is
if there exists an ergodic measure of full dimension. In a restricted domain
of parameter values we found generically an absolutely continuous ergodic measure which obviously has dimension two.
\begin{theorem}
For almost all $(\beta_{1},\beta_{2})\in (0,0.649)$ with $\beta_{1}+\beta_{2}\ge 1$ and $\beta_{1}\beta_{2}\ge 1/4$ there is an absolutely continuous ergodic measure
for $([-1,1]^2,f_{\beta_{1},\beta_{2}})$.
\end{theorem}
This theorem mainly is a consequence of our results about
overlapping self-similar measures one
the real line \cite{[NE1]}. We will
construct the measure of full dimension as a product of an overlapping
self-similar with normalized Lebesgue measure. From \cite{[NE1]} we
then deduce absolute continuity of this measure. We do not know if the condition $\beta_{1},\beta_{2}<0.649$ in theorem 2.1 is necessary, in fact
it is due to the techniques we used in \cite{[NE1]}.
On the other hand from our second theorem we see that the condition $\beta_{1}\beta_{2}\ge 1/4$ in
theorem 2.1 is necessary.
\begin{theorem}
For $(\beta_{1},\beta_{2})\in(0,1)$
with $\beta_{1}+\beta_{2}\ge 1$ and $\beta_{1}\beta_{2}< 1/4$
we have
\[ \sup\{\dim_{H}\mu|\mu~~f_{\beta_{1},\beta_{2}}\mbox{-ergodic}\}<2. \]
\end{theorem}
This example shows that it is not always possible to
find the Hausdorff dimension of an invariant set by constructing an ergodic  measure
of full Hausdorff dimension.
Roughly speaking the reason why there is not always an ergodic measure of full Hausdorff dimension here
is that one can not maximize the stable and the unstable dimension (the dimension of
conditional measures on partitions in stable resp. unstable directions) at the same time.
In another context this praenomen was observed before by Manning and McClusky \cite{[MM]}.
\\\\
Now consider the for a moment the Fat Baker's transformation $f_{\beta}:=f_{\beta,\beta}$ with
$\beta\in(0.5,1)$. It follows from the work of Alexander and Yorke \cite{[AY]} together with Solymak's theorem on Bernoulli convolutions \cite{[SO1]}
that for almost all $\beta\in(0.5,1)$ we have $\dim_{H}\mu_{SRB}=2$ where $\mu_{SRB}$ is
the Sinai-Ruelle-Bowen measure for the system $([-1,1]^{2},f_{\beta})$, see \cite{[NE2]}.
This means that in the symmetric situation, in contrast to the asymmetric case, we generically have an ergodic measure of full dimension in the whole parameter domain.
\section{Symbolic coding and representation of ergodic measures}
Let $\Sigma=\{-1,1\}^{\Z}$, $\Sigma^{+}=\{-1,1\}^{\N_{0}}$
and $\Sigma^{-1}=\{-1,1\}^{\Z^{-}}$. The forward shift map $\sigma$ on $\Sigma$ (resp. $\Sigma^{+}$) is given by
$\sigma((s_{k}))=(s_{k+1})$ and the system $(\Sigma,\sigma)$
(resp. $(\Sigma^{+},\sigma)$) is know as full shift on two symbols \cite{[KH]}.
Given $\underline{s}\in\Sigma^{+}$ we denote by $\sharp_{k}(\underline{s})$ the cardinality of $\{i|s_{i}=-1~,~i=0\dots k\}$. \\\\
For $\beta_{1},\beta_{2}\in (0,1)$ with $\beta_{1}+\beta_{2}\ge 1/2$ we now define a map $\pi_{\beta_{1},\beta_{2}}$ from $\Sigma^{+}$
onto $[-1,1]$ in the following way.
Let
\[  \pi^{\ast}_{\beta_{1},\beta_{2}}(\underline{s})=\sum_{k=0}^{\infty}s_{k}\beta_{2}^{\sharp_{k}(\underline{s})}\beta_{1}^{k-\sharp_{k}(\underline{s})+1}
.\]
We scale this map so that it is onto $[-1,1]$.
by be the affine transformation  $L_{\beta_{1},\beta_{2}}$ on the line that maps $\frac{-\beta_{2}}{1-\beta_{2}}$ to $-1$ and
$\frac{\beta_{1}}{1-\beta_{1}}$ to $1$;
$\pi_{\beta_{1},\beta_{2}}
=L_{\beta_{1},\beta_{2}}\circ \pi^{\ast}_{\beta_{1},\beta_{2}}$. \\
Now define the maps $\varsigma$ from $\Sigma^{-}$ onto
[-1,1]  corresponding to the signed dyadic expansion of a number
by
\[ \varsigma(\underline{s} ) = \sum_{k=1}^{\infty}s_{-k}2^{-k}
~\mbox{ where }~\underline{s}=(s_{k})_{k\in\Z^{-}}\in\Sigma^{-}.\]
We are now able to define the coding map
for the systems $([-1,1]^{2},f_{\beta_{1},\beta_{2}})$ by
\[ \bar\pi_{\beta_{1},\beta_{2}}:\Sigma \longmapsto
[-1,1]^{2} \qquad\mbox{with} \qquad \pi_{\beta_{1},\beta_{2}}((s_{k}))
=(\pi_{\beta_{1},\beta_{2}}((s_{k})_{k\in\N_{0}}),
\varsigma((s_{k})_{k\in\Z^{-}})).\]
Obviously $\bar\pi_{\beta_{1},\beta_{2}}$ is onto and continuous if we endow
$\Sigma$ with the natural product topology. Moreover we have
\begin{proposition}
$\bar\pi_{\beta_{1},\beta_{2}}$
conjugates the backward shift $\sigma^{-1}$ and $f_{\beta_{1},\beta_{2}}$ i.e.
\[ f_{\beta_{1},\beta_{2}}\circ\pi_{\beta_{1},\beta_{2}}=\pi_{\beta_{1},\beta_{2}}\circ \sigma^{-1}\]
on
\[ \bar\Sigma=(\Sigma\backslash\{(s_{k})|
\exists k_{0}
\forall k \le k_{0} :s_{k}=1 \rbrace)\cup\{(1)\}.\]
\end{proposition}
\proof
Let $\underline{s}=(s_{k}) \in \bar\Sigma$. We have $(s_{k+1})_{k \in \bf{Z}^{-}}\not=(\dots,1,1,-1)$ and hence
\[ \varsigma((s_{k+1})_{k\in\Z^{-}})=
\sum_{k=1}^{\infty}s_{-k+1}2^{-k} \ge 0 \Leftrightarrow s_{0}=1.\]
Thus
\[ f_{\beta_{1},\beta_{2}} \circ \bar\pi_{\beta_{1},\beta_{2}}((s_{k+1}))=\]
\[ \lbrace
\begin{array}{cc}
( \beta_{1}\pi_{\beta_{1},\beta_{2}}((s_{k+1})_{k\in\N_{0}})+(1-\beta_{1}),
2\varsigma((s_{k+1})_{k\in\Z^{-}})-1)
\mbox{ if } s_{0}=1 \\
( \beta_{2}\pi_{\beta_{1},\beta_{2}}((s_{k+1})_{k\in\N_{0}})-(1-\beta_{2}),
2\varsigma((s_{k+1})_{k\in\Z^{-}})+1)
\mbox{ if } s_{0}=-1
\end{array}.\]
On the other hand we have
\[ \pi_{\beta_{1},\beta_{2}}(\sigma(\underline{s}))
=L_{\beta_{1},\beta_{2}}(\pi^{\ast}_{\beta_{1},\beta_{2}}((s_{k+1})))={\lbrace
{\begin{array}{cc}
 L_{\beta_{1},\beta_{2}}(\beta_{1}^{-1}\pi^{\ast}_{\beta_{1},\beta_{2}}((s_{k}))-1) \quad \mbox{if} \quad s_{0}=1   \\
~~ L_{\beta_{1},\beta_{2}}(\beta_{2}^{-1}\pi^{\ast}_{\beta_{1},\beta_{2}}((s_{k}))+1)
\quad \mbox{if} \quad s_{0}=-1
\end{array} }}\]
\[={\lbrace
{\begin{array}{cc}
 \beta_{1}^{-1}\pi_{\beta_{1},\beta_{2}}(\underline{s})+(1-\beta_{1}^{-1}) \quad \mbox{if} \quad s_{0}=1   \\
~~\beta_{2}^{-1}\pi_{\beta_{1},\beta_{2}}(\underline{s})-(1-\beta_{2}^{-1})
\quad \mbox{if} \quad s_{0}=-1
\end{array} }} .\]
By these equations and the definition of $\varsigma$ we now see that
$f_{\beta_{1},\beta_{2}} \circ \bar\pi_{\beta_{1},\beta_{2}}((s_{k+1}))=\bar\pi_{\beta_{1},\beta_{2}}((s_{k}))$. $\sigma$ as a map of $\Sigma$ is invertible and we get $f_{\beta_{1},\beta_{2}} \circ \bar\pi_{\beta_{1},\beta_{2}}(\underline{s})=
\bar\pi_{\beta_{1},\beta_{2}}(\sigma^{-1}(\underline{s}))$ for all $\underline{s}\in\bar\Sigma$.
\begin{flushright}
$\Box$
\end{flushright}
Using our symbolic coding we can describe all ergodic measures for $f_{\beta_{1},\beta_{2}}$. To this end we introduce the following
notation: $M(X,f)$ denotes the space of all $f$-ergodic Borel probability measures
on $X$. It is well known in ergodic theory that if $X$ is compact
$M(X,f)$ is a nonempty convex $weak^{*}$ compact metricable space, \cite{[WA]} or \cite{[DGS]}.
\begin{proposition}
$\mu\longmapsto\bar\mu_{\beta_{1},\beta_{2}}:=\mu\circ\bar\pi_{\beta_{1},\beta_{2}}^{-1}$ is a continuous affine map from $M(\Sigma,\sigma)$
onto $M([-1,1]^{2},f_{\beta_{1},\beta_{2}})$.
\end{proposition}
\proof
It is obvious by standard arguments in measure theory \cite{[MT]} that
the map in question is continuous and affine since $\bar\pi_{\beta_{1},\beta_{2}}$ is continuous.
If $\mu$ is shift ergodic
we have $\mu(\bar\Sigma)=1$. We know
from Proposition 3.1 that
$\bar\pi_{\beta_{1},\beta_{2}}$ conjugates the backward shift
and $f_{\beta_{1},\beta_{2}}$ on $\bar\Sigma$ hence we get
that $\bar\mu_{\beta_{1},\beta_{2}}$ is $f_{\beta_{1},\beta_{2}}$-ergodic.
It remains to show that the map is
onto
$M([-1,1]^{2},f_{\beta_{1},\beta_{2}})$. This is a not completely trivial
exercise in ergodic theory. Let us choose an arbitrary measure $\xi$  in $M([-1,1]^{2},f_{\beta_{1},\beta_{2}})$.\newline
We first want to show that $\xi(\pi_{\beta_{1},\beta_{2}}(\Sigma\backslash\bar\Sigma))=0$.
Let $D$ be set of all numbers of the form $k/2^{n}$ with
$n\in \N$ and $|k|\le n-1$.
A direct calculation shows that
\[ \pi_{\beta_{1},\beta_{2}}(\Sigma\backslash\bar\Sigma)
=(D\times[-1,1])\cup(\{1\} \times [-1,1))
=(\bigcup_{k=0}^{\infty}f_{\beta_{1},\beta_{2}}^{-k}(\{0\}
\times [1,-1]))\cup(\{1\} \times [-1,1))
.\]
Recall that the measure $\xi$ is in particular shift invariant. Hence the measure of the first set in union is zero because it is given by a disjunct infinite union of sets with the same measure. The measure of
the second set is zero since $\{1\}\times[-1,1)\subseteq
f_{\beta_{1},\beta_{2}}^{-k}(\{1\}\times[1-2\beta_{1}^{k},1))~~ \forall~k\ge 0$. \newline
Now take a Borel probability measure $\mu_{pre}$
such that $ \mu_{pre}\circ\pi_{\beta_{1},\beta_{2}}^{-1}            =\xi$.
$\mu_{pre}$ is
not necessary shift invariant so we define a measure
$\mu$ as a $\mbox{weak}^{*}$ accumulation point of the sequence
\[\mu_{n}:=\frac{1}{n+1}\sum_{i=0}^{n}\mu_{pre}\circ\sigma^{-n}.\]
From the considerations above we have $\mu_{pre}(\bar\Sigma)=1$ and hence:
 \[ \mu_{n}\circ\pi_{\beta_{1},\beta_{2}}^{-1}
 =\frac{1}{n+1}\sum_{i=0}^{n}\mu_{pre}\circ\sigma^{-n}\circ
 \pi_{\beta_{1},\beta_{2}}^{-1}\]
\[ =\frac{1}{n+1}\sum_{i=0}^{n}\mu_{pre}\circ
 \pi_{\beta_{1},\beta_{2}}^{-1}\circ f_{\beta_{1},\beta_{2}}^{-i}
=\frac{1}{n+1}\sum_{i=0}^{n}\xi \circ f_{\beta_{1},\beta_{2}}^{-i}
=\xi. \]
Thus $\bar\mu_{\beta_{1},\beta_{2}}$ is just the measure $\xi$
and $\mu$ is shift invariant by definition. We have thus shown
that the set $M(\xi):=\{\mu|\mu$ $\sigma$-invariant and $\mu_{\beta_{1},\beta_{2}}=\xi\}$
of Borel measures on $\Sigma$ is not empty. Since the map $\mu\longmapsto\bar\mu_{\beta_{1},\beta_{2}}$ is continuous and affine on the
set of $\sigma$-invariant measures we know that $M(\xi)$ is compact and convex. It is a consequence of Krein-Milman theorem that there
exists an extremal point $\mu$ of $M(\xi)$.\newline We claim that $\mu$ is an extremal point of the set of all $\sigma$-invariant Borel measures on $\Sigma$ and hence ergodic.\newline If this is not the case then we have $\mu=t\mu_{1}+(1-t)\mu_{2}$ where $t\in(0,1)$ and $\mu_{1},\mu_{2}$ are two distinct $\sigma$-invariant measures. This implies
$\xi=t(\mu_{1})_{\beta_{1},\beta_{2}}+(1-t)(\mu_{2})_{\beta_{1},\beta_{2}}$. Since $\xi$ is ergodic we have $(\mu_{1})_{\beta_{1},\beta_{2}}=(\mu_{2})_{\beta_{1},\beta_{2}}=\xi$
and hence $\mu_{1},\mu_{2}\in M(\xi)$. This is a contradiction to
$\mu$ being extremal in $M(\xi)$.
\begin{flushright}
$\Box$
\end{flushright}
\section{Construction of absolutely continuous ergodic measures}
We now construct absolutely continuous ergodic measures for the systems $([-1,1]^{2},f_{\beta_{1},\beta_{2}})$. Let $b$ denote the Bernoulli measure on
the shift $\Sigma$ (resp. $\Sigma^{+}$ or $\Sigma^{-}$, which is the product of the discrete
measure giving $1$ and $-1$ the probability $1/2$.
The Bernoulli measure is ergodic with respect to forward and
backward shifts, see \cite{[KH]}.\newline
Given $b$ on
$\Sigma^{-}$
we set $\ell^{p}=b\circ\varsigma^{-1}$. $\ell$
is the normalized Lebesgue measure on [-1,1].\newline
Given $b$ on $\Sigma^{+}$
we define two Borel probability
measures on the real line
by \[b^{\ast }_{\beta_{1},\beta_{2}}=b\circ(\pi^{\ast} _{\beta_{1},\beta_{2}})^{-1}
~~~\mbox{and}~~~
b_{\beta_{1},\beta_{2}}=b\circ(\pi _{\beta_{1},\beta_{2}})^{-1}.
\]
The measure $b_{\beta_{1},\beta_{2}}$ is just
$b^{\ast}_{\beta_{1},\beta_{2}}$ scaled on the interval
$[-1,1]$ by the transformation $L_{\beta_{1},\beta_{2}}$: \[ b_{\beta_{1},\beta_{2}}=b^{\ast}_{\beta_{1},\beta_{2}}\circ L_{\beta_{1},\beta_{2}}^{-1}. \]
In the following proposition we describe an ergodic measure fore the generalized Bakers transformations using
the Bernoulli measure $b$.
\begin{proposition}
For all $\beta_{1},\beta_{2}\in(0,1)$ with $\beta_{1}+\beta_{2}\ge 1$ we have
\[ \bar b_{\beta_{1},\beta_{2}}:=b\circ\bar \pi_{\beta_{1},\beta_{2}}^{-1}=b_{\beta_{1},\beta_{2}}\times\ell \in M([-1,1]^{2},f_{\beta_{1},\beta_{2}}). \]
\end{proposition}
\proof By Proposition 3.2  we get that $b\circ\bar \pi_{\beta_{1},\beta_{2}}^{-1}$ is $f_{\beta_{1},\beta_{2}}$-ergodic since
$b$ is $\sigma-$ergodic. Moreover by the product structure of $\bar \pi_{\beta_{1},\beta}$ we have
\[ b\circ\bar \pi_{\beta_{1},\beta_{2}}^{-1}=b\circ \pi_{\beta_{1},\beta_{2}}\times b\circ\iota^{-1}=b_{\beta_{1},\beta_{2}}\times \ell     \]
where we use the fact that the Bernoulli measure $b$ on $\Sigma$ is a product of $b$ on $\Sigma^{+}$ with $b$ on
$\Sigma^{-1}$
\begin{flushright}
$\Box$
\end{flushright}
The measures
$ b^{\ast}_{\beta_{1},\beta_{2}}$ are by definition a special class of overlapping self similar measures studied in \cite{[NE1]}. The following proposition is
just a consequence of Theorem I of \cite{[NE1]}
\begin{proposition}
For almost all $(\beta_{1},\beta_{2})\in (0,0.649)$ with $\beta_{1}+\beta_{2}\ge 1$ and $\beta_{1}\beta_{2}\ge 1/4$ the
measure $b^{\ast}_{\beta_{1},\beta_{2}}$ is absolutely continuous.
\end{proposition}
Now the proof of Theorem 2.2 is obviously.\\\\
{\bf Proof of Theorem 2.1.}~
In the relevant parameter domain we
generically have absolute continuity of $b^{*}_{\beta_{1},\beta_{2}}$ by Proposition 4.2.
Since  $b_{\beta_{1},\beta_{2}}= b^{*}_{\beta_{1},\beta_{2}}\circ L_{\beta_{1},\beta_{2}}^{-1}$ this clearly implies  absolute continuity of $b_{\beta_{1},\beta_{2}}$. Now
Proposition 4.1 implies absolute continuity of the measure $\bar b_{\beta_{1},\beta_{2}}$ which is $f_{\beta_{1},\beta_{2}}$-ergodic.
\begin{flushright}
$\Box$
\end{flushright}

\section{Dimension estimates on all ergodic measures}
In this section we proof two upper bounds on the dimension
of all $f_{\beta_{1},\beta_{2}}$-ergodic measures using metric entropy of these measures. Theorem 2.2 will be a consequence.\\
Given $\mu \in M(\Sigma,\sigma)$ we denote by $h(\mu)$ the metric entropy of $\mu$. We refer to \cite{[WA]} or \cite{[KH]} the definition
and the properties of this quantity.
\begin{proposition}
For all $\mu \in M(\Sigma,\sigma)$ and all $\beta_{1},\beta_{2}\in(0,1)$ with $\beta_{1}+\beta_{2}\ge 1$ we have
\[\dim_{H}\bar \mu_{\beta_{1},\beta_{2}}\le \frac{h(\mu)}{\log 2}+1\]
\end{proposition}
\proof
Let $\tilde\mu_{\beta_{1},\beta_{2}}$ by the projection of the measure  $\bar\mu_{\beta_{1},\beta_{2}}$ onto
the second coordinate axis. Since $\dim_{H}(B\times[-1,1])=\dim_{H}B+1$ for all sets $B$ we have
\[\dim_{H} \bar\mu_{\beta_{1},\beta_{2}}\le\dim \tilde\mu_{\beta_{1},\beta_{2}}+1. \]
just by the definition of the Hausdorff dimension of a measure. Now we have to estimate the dimension of the projection.
By
$\tilde\mu_{\beta_{1},\beta_{2}}=\bar\mu_{\beta_{1},\beta_{2}}\circ pr_{y}^{-1}=\mu\circ \bar\pi_{\beta_{1},\beta_{2}}^{-1}\circ pr_{y}^{-1}$
and the product structure of the coding map $\bar\pi_{\beta_{1},\beta_{2}}$ we see that
$\tilde\mu_{\beta_{1},\beta_{2}}$ is ergodic with respect to the map
$$
{ f(y) } = \lbrace
\begin{array}{cc}
2y-1\quad \mbox{if}\quad y \geq 0   \\
2y+1\quad \mbox{if} \quad y<0.
\end{array}
$$
Thus the Hausdorff dimension of $\tilde\mu_{\beta_{1},\beta_{2}}$ is well known (see \cite{[PE]})
$$ \dim_{H}\tilde\mu_{\beta_{1},\beta_{2}}=\frac{h(\tilde\mu_{\beta_{1},\beta_{2}})}{\log 2}. $$
Moreover we know that $([-1,1],f,\tilde\mu_{\beta_{1},\beta_{2}})$ is a measure theoretical factor of $([-1,1]^{2},$ $f_{\beta_{1},\beta_{2}},\bar\mu_{\beta_{1},\beta_{2}})$.
and that this system is a factor of $(\Sigma,\sigma,\mu)$. Hence we get by well known properties of the entropy (see \cite{[DGS]}) $h(\tilde\mu_{\beta_{1},\beta_{2}})\le h(\bar\mu_{\beta_{1},\beta_{2}})\le h(\mu)  $
which completes the proof.

\begin{flushright}
$\Box$
\end{flushright}
To state the other estimate we need a few notation. Let $pr$ we the projection from $\Sigma$ to $\Sigma^{+}$.
Given $\mu$ in $M(\Sigma,\sigma)$ we define $\hat\mu \in M(\Sigma^{+},\sigma)$ by $\hat\mu=\mu\circ pr^{-1}$.
Moreover set \[\Xi_{\beta_{1},\beta_{2}}(\hat\mu)=-(\hat\mu(\{\underline{s}\in\Sigma^{+}|s_{0}=1\})\log\beta_{1}+\hat\mu(\{\underline{s}\in\Sigma^{+}|s_{0}=-1\})\log\beta_{2}\]
With these notations we have
\begin{proposition}
For all $\mu \in M(\Sigma,\sigma)$ and all $\beta_{1},\beta_{2}\in(0,1)$  with $\beta_{1}+\beta_{2}\ge 1$ we have
\[\dim_{H}\bar \mu_{\beta_{1},\beta_{2}}\le \frac{h(\hat\mu)}{\Xi_{\beta_{1},\beta_{2}}(\hat\mu)}+1 \]

\end{proposition}
\proof
The proof of this proposition has several steps. First we show the following inequality
\[ \dim_{H}\bar\mu_{\beta_{1},\beta_{2}}\le \dim_{H}\hat\mu_{\beta_{1},\beta_{2}}+1 \]
Let $B$ be an arbitrary Borel set with $\hat\mu_{\beta_{1},\beta_{2}}(B)=1$.
Since the projection of $\bar\mu_{\beta_{1},\beta_{2}}$ onto the first
coordinate axis is $\hat\mu_{\beta_{1},\beta_{2}}$ we get
$\bar\mu_{\beta_{1},\beta_{2}}(B\times[-1,1])=1$ Thus
\[ \dim_{H}\bar\mu_{\beta_{1},\beta_{2}}\le\dim_{H}(B\times[-1,1])=\dim_{H}(B)+1\]
Now our claim follows just by the definition of the Hausdorff dimension of a measure.\\\\
Now we have to estimate the dimension of
the projected measure; we have to show that
\[\dim_{H}\hat \mu_{\beta_{1},\beta_{2}}\le \frac{h(\hat\mu)}{\Xi_{\beta_{1},\beta_{2}}(\hat\mu)} \]
Let us define a metric $\delta^{\beta_{1},\beta_{2}}$ on
$\Sigma^{+}$ by
\[ \delta^{\beta_{1},\beta_{2}}(\underline{s},\underline{t})
=\beta_{1}^{|\underline{s}\wedge \underline{t}|-\sharp_{|\underline{s}\wedge \underline{t}|-1}(\underline{s})}
\beta_{2}^{\sharp_{|\underline{s}\wedge \underline{t}|-1}(\underline{s})}.\]
where  $\sharp_{k}(\underline{s})$ is the cardinality of $\{i|s_{i}=-1~,i=0\dots k\}$ and $|\underline{s}\wedge \underline{t}\|=\min\{i|s_{i}\not=t_{i}\}$.
Now we claim that
\[ d^{\beta_{1},\beta_{2}}(\underline{s},\hat\mu):=\lim_{\epsilon\longrightarrow 0}\frac{\log\hat\mu( B_{\epsilon}^{\beta_{1},\beta_{2}})(\underline{s})}{\log\epsilon}=\frac{h_{\hat\mu}(\sigma)}{\Xi_{\beta_{1},\beta_{2}}(\hat\mu)}\qquad\hat\mu\mbox{-almost everywhere}.\]
Here $d^{\beta_{1},\beta_{2}}$ is the local dimension of
the measure $\hat\mu$ with
respect to metric $\delta^{\beta_{1},\beta_{2}}$ and accordingly
$B_{\epsilon}^{\beta_{1},\beta_{2}}$ is a ball of radius
$\epsilon$ with respect to this metric. Applying Birkhoff's
ergodic theorem (see 4.1.2. of \cite{[KH]}) to $(\Sigma^{+},\sigma,\hat\mu)$
with the function
\[h(\underline{s})=
{\lbrace
{\begin{array}{cc}
  \log \beta_{1}\quad\mbox{if} \quad s_{0}=1   \\
  ~~\log \beta_{2}\quad\mbox{if} \quad s_{0}=-1
\end{array} }}
\]
we see that\[ \lim_{n\longrightarrow\infty}\frac{1}{n+1}\log \mbox{diam}_{\beta_{1},\beta_{2}}([s_{0},\dots,s_{n}])= \lim_{n\longrightarrow\infty}\frac{1}{n+1}\sum_{k=0}^{n+1}h(\sigma^{k}(\underline{s}))=\int h~d\hat\mu (\underline{s})=\Xi_{\beta_{1},\beta_{2}}(\hat\mu)\]
$\hat\mu$-almost everywhere. By Shannon-McMillan-Breiman theorem (see \cite{[DGS]} 13.4.) we have:
\[\lim_{n\longrightarrow\infty}-\frac{1}{n+1}\log\hat\mu([s_{0},\dots,s_{n}])=h_{\hat\mu}(\sigma)\qquad\hat\mu\mbox{-almost everywhere}.\]
Thus we see:
\[\lim_{\epsilon\longrightarrow 0}\frac{\log B_{\epsilon}^{\beta_{1},\beta_{2}}(\underline{s})}{\log\epsilon} =\lim_{n\longrightarrow\infty}
\frac{\log\hat\mu([s_{0},\dots,s_{n}])}{\mbox{diam}_{\beta_{1},\beta_{2}}([s_{0},\dots,s_{n}]_{0}) }=\frac{h(\hat\mu)}{\Xi_{\beta_{1},\beta_{2}}(\hat\mu)} \]
By Billigsly's Lemma about the relation of local and global dimension \cite{[PE]} we get
\[ \dim_{H}^{\beta_{1},\beta_{2}}\hat\mu=\frac{h(\hat\mu)}{-\Xi_{\beta_{1},\beta_{2}}(\hat\mu)} \]
where the Hausdorff dimension $\dim_{H}^{\beta_{1},\beta_{2}}$ has to be calculated with respect to the metric  $\delta^{\beta_{1},\beta_{2}}$. \newline\newline
Now we claim that the map $\pi^{*}_{\beta_{1},\beta_{2}}$ is Lipschitz with respect to the metric $\delta^{\beta_{1},\beta_{2}}$
\[|\pi^{*}_{\beta_{1},\beta_{2}}(\underline{s})-\pi^{*}_{\beta_{1},\beta_{2}}(\underline{t})|\le
\sum_{k=|\underline{s}\wedge \underline{t}|}^{\infty}
|s_{k}\beta_{1}^{k-\sharp_{k}(\underline{s})+1}\beta_{2}^{\sharp_{k}(\underline{s})}-
t_{k}\beta_{1}^{k-\sharp_{k}(\underline{t})+1}\beta_{2}^{\sharp_{k}(\underline{t})}|
\]
\[
=\beta_{1}^{|\underline{s}\wedge \underline{t}|-\sharp_{|\underline{s}\wedge \underline{t}|-1}(\underline{s})}
\beta_{2}^{\sharp_{|\underline{s}\wedge \underline{t}|-1}(\underline{s})}\]
\[
\sum_{k=0}^{\infty}
|s_{k+|\underline{s}\wedge \underline{t}|}\beta_{1}^{k-\sharp_{k}(\sigma^{|\underline{s}\wedge-\underline{t}|}(\underline{s}))+1}\beta_{2}^{\sharp_{k}(\sigma^{|\underline{s}\wedge\underline{t}|}(\underline{s}))}-
t_{k+|\underline{s}\wedge \underline{t}|}\beta_{1}^{k-\sharp_{k}(\sigma^{|\underline{s}\wedge\underline{t}|}(\underline{t}))+1}\beta_{2}^{\sharp_{k}(\sigma^{|\underline{s}\wedge\underline{t}|}(\underline{t}))}|
\]
\[\le \delta^{\beta_{1},\beta_{2}}(\underline{s},\underline{t}) \frac{2}{1-\max\{\beta_{1},\beta_{2}\}}.\]
The the map $\pi_{\beta_{1},\beta_{2}}$ is just $\pi^{*}_{\beta_{1},\beta_{2}}$ scaled on $[-1,1]$ and hence
as Lipschitz with respect $\delta^{\beta_{1},\beta_{2}}$.
Recall that $\hat\mu_{\beta_{1},\beta_{2}}=\hat\mu\circ \pi_{\beta_{1},\beta_{2}}$.
Since applying a Lipschitz map
to the measures $\hat\mu$ does obvious not increase its Hausdorff dimension, the proof is complete.
\begin{flushright}
$\Box$
\end{flushright}
Combining proposition 5.1 and 5.2 we now proof Theorem 2.3\\\\
{\bf Proof of Theorem 2.3.}\quad
If $\beta_{1}\beta_{2}< 0.25$ then we have $h(b)<\Xi_{\beta_{1},\beta_{2}}(b)$. By upper semi continuity of the metric entropy there is
a $\mbox{weak}^{*}$ neighborhood $U$ of $b$ in $M(\Sigma^{+},\sigma)$ such that $h(\mu)/\Xi_{\beta_{1},\beta_{2}}(\mu)\le c_{1}<1$
holds for all $\mu\in U$.
By Proposition 5.1 we get
$\dim_{H}\bar\mu_{\beta_{1},\beta_{2}}\le c+1<2$ for all
$\mu \in \tilde U=pr^{-1}(U)$. Obviously $\tilde U$ is
a neighborhood of $b$ in $M(\Sigma,\sigma)$. Furthermore
by Proposition 5.1 and upper-semi continuity of metric entropy we get that
$\dim_{H}\bar\mu_{\beta_{1},\beta_{2}}\le c_{2}+1<2$ for all
$\mu\in M(\Sigma,\sigma)\backslash \bar U$. Putting these facts together we get
\[ \dim_{H}\bar\mu_{\beta_{1},\beta_{2}}\le\max\{c_{1},c_{2}\}+1<2=\dim[-1,1]^{2} \qquad \forall \mu\in M(\Sigma,\sigma).\]
\begin{flushright}
$\Box$
\end{flushright}

\end{document}